\documentclass[11pt]{article}

\usepackage{bridges}
\usepackage{amsfonts,amssymb,amsthm,eucal,amsmath}
\usepackage{graphicx}
\usepackage[T1]{fontenc}
\usepackage{latexsym,url}
\usepackage{enumerate}
\usepackage{subfig}
\usepackage{wrapfig}
\usepackage{pinlabel}
\usepackage{color}
\usepackage{microtype}
\usepackage[
hidelinks,         
pdftex]{hyperref}

\newtheorem*{no_num_thm}{Theorem}

\theoremstyle{definition}

\theoremstyle{remark}

\newcommand{\reffig}[1]{Figure~\ref{Fig:#1}}
\newcommand{\refsec}[1]{Section~\ref{Sec:#1}}

\newcommand{\CC}{\mathbb{C}}

\newcommand{\RR}{\mathbb{R}}

\title{Triple gear\thanks{This work is in the public domain.}}
\renewcommand\footnotemark{}

\author{
\begin{tabular}{cc}
Saul Schleimer            & Henry Segerman \\
Mathematics Institute     & Department of Mathematics \\
University of Warwick     & Oklahoma State University \\
Coventry CV4 7AL          & Stillwater, OK 74078 \\
United Kingdom            & USA \\
s.schleimer@warwick.ac.uk & henry@segerman.org
\end{tabular}}

\date{}

\begin{document}
\maketitle
\begin{abstract}
A relatively common sight in graphic designs is a planar arrangement
of three gears in contact.  However, since neighboring gears must
rotate in opposite directions, none of the gears can move.  We give a
non-planar, and non-frozen, arrangement of three linked gears.
\end{abstract}


\section{Introduction}

\begin{wrapfigure}[11]{l}{0.3\textwidth}
\centering
\vspace{-10pt}
\labellist
\footnotesize\hair 2pt
\pinlabel \textcolor{red}{meridian} [Br] at 185 407 
\pinlabel \textcolor{red}{longitude} [tr] at 10 190
\definecolor{darkgreen}{rgb}{0,0.659,0}
\pinlabel \textcolor{darkgreen}{$(1,1)$ curve} at 310 160 
\endlabellist
\includegraphics[width=0.3\textwidth]{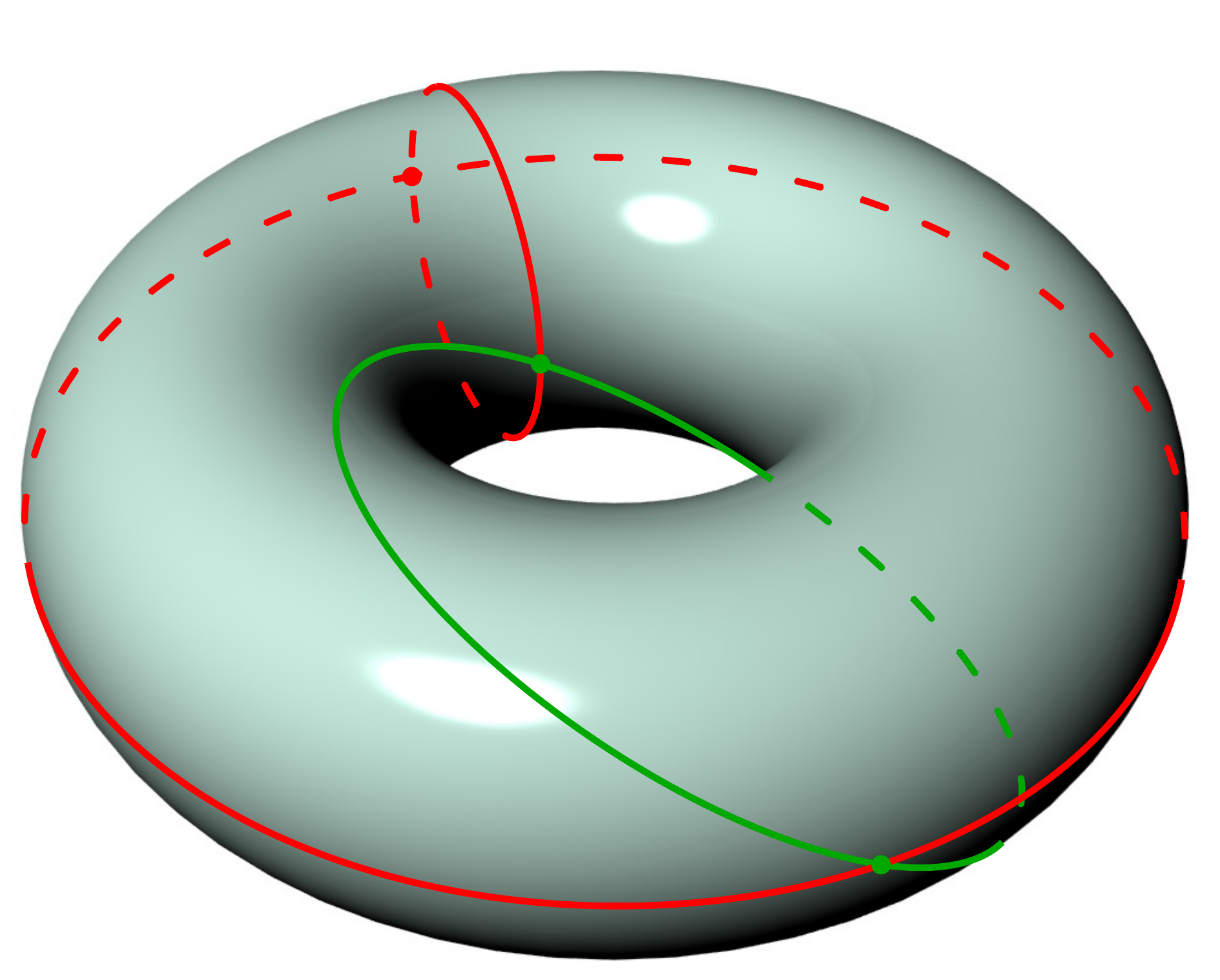}
\caption{Curves on a torus.} 
\label{Fig:Torus curves}
\end{wrapfigure}

The mathematician's \emph{torus} is simply the crust of a bagel: a
surface with one hole.  There are two special kinds of curves on the
torus called \emph{meridians} and \emph{longitudes}.  These are the
curves obtained by slicing the bagel vertically or horizontally,
respectively.  We call a curve on the torus a \emph{$(p,q)$ curve} if
it crosses a longitude $p$ times and a meridian $q$ times.  See
\reffig{Torus curves}.

\begin{wrapfigure}[12]{r}{0.32\textwidth}
\centering
\vspace{-10pt}
\includegraphics[width=0.32\textwidth]{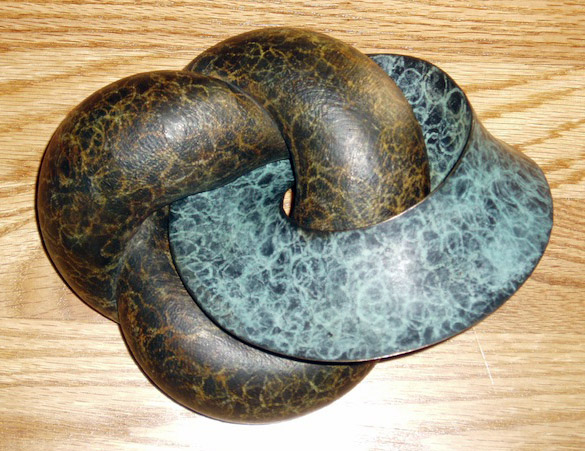}
\caption{\emph{Umbilic Rolling Link}, by Helaman Ferguson.}
\label{Fig:Rolling}
\end{wrapfigure}

Helaman Ferguson's \emph{Umbilic Torus NC} is a classic work of
mathematical art; topologically it is a torus with a ridge along a
$(1,3)$ curve~\cite{Ferguson}.  He has also designed a complementary
sculpture, called the \emph{Umbilic Torus NIST}; here we find a torus
with a valley along a $(3,1)$ curve.  These complementary tori are
mated in his piece \emph{Umbilic Rolling Link}, see \reffig{Rolling}.
Note that the ridge of the NC torus meshes with the valley of the NIST
torus.  Thus if one of the tori rotates along itself then so must the
other, with speeds in a ratio of $3:1$.

\begin{wrapfigure}[13]{l}{0.32\textwidth}
\centering
\vspace{-7pt}
\includegraphics[width=0.32\textwidth]{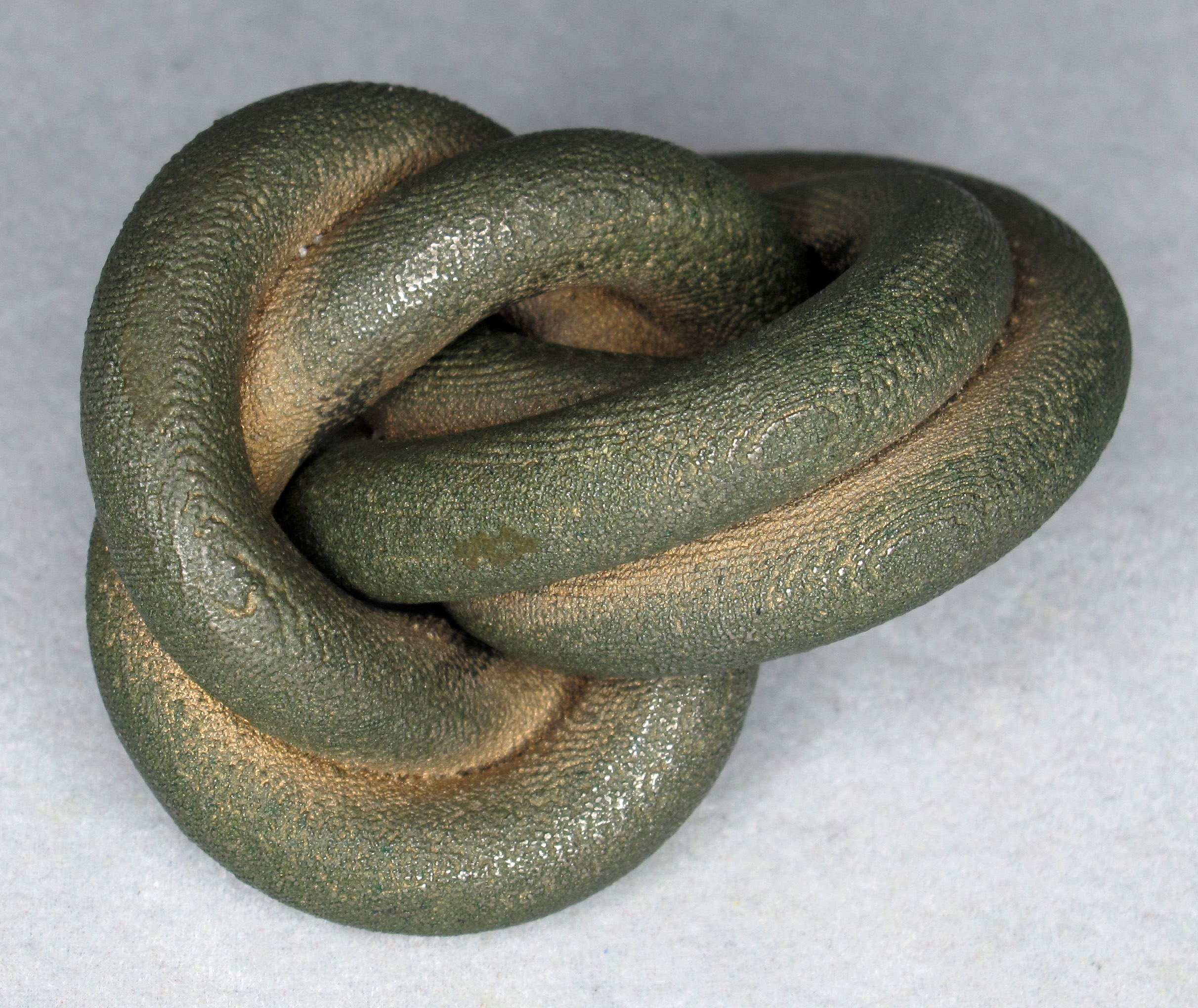}
\caption{\emph{Knotted Gear}, by Oskar van Deventer.} 
\label{Fig:Knotted_Gear}
\end{wrapfigure}

A closely related piece by Oskar van Deventer inspired this project.
His \emph{Knotted Gear}\footnote{See
  \url{http://www.youtube.com/watch?v=ZqiHHlE1SD8} for a video of the
  knotted gear in motion.} consist of two rings, a thickened $(-3,2)$
curve (a left-handed trefoil) linked with a thickened $(2,-3)$ curve.
See \reffig{Knotted_Gear}.  The two rings mesh with speeds in a ratio
of $3:2$.  We quickly decided that by varying $p$ and $q$ any rational
ratio could, in principle, be achieved.  We next asked ``Can this be
done with three or more rings?''

To formalize this problem we begin with some examples and definitions.
An example of a rigid body is a wheel.  An example of a rigid motion
is a rotation of the wheel about its center through some angle.  A
\emph{movement} is a one-parameter family of rigid motions $\{ g_t
\}_{t \in \RR}$.  For example, consider spinning the wheel about its
center with constant angular velocity $\omega$: then at time $t$, the
motion $g_t$ is a rotation through $\omega t$ degrees. Given two
motions $m_1$ and $m_2$, the \emph{composition of $m_1$ with $m_2$},
denoted $m_2 \circ m_1$, is the motion we get by performing $m_1$
followed by $m_2$.  A movement $g$ is \emph{simple} if for any times
$s$ and $t$ we have $g_{s+t} = g_s \circ g_t$.  Note that the simple
movements are either stationary, translational, rotational or screw
movement.  These movements get their name from the classical simple
machines: the lever, the wheel, the screw, and so on.  For example,
the spinning wheel is a simple movement.

\begin{wrapfigure}[11]{r}{0.32\textwidth}
\vspace{-20pt}
\includegraphics[width=0.32\textwidth]{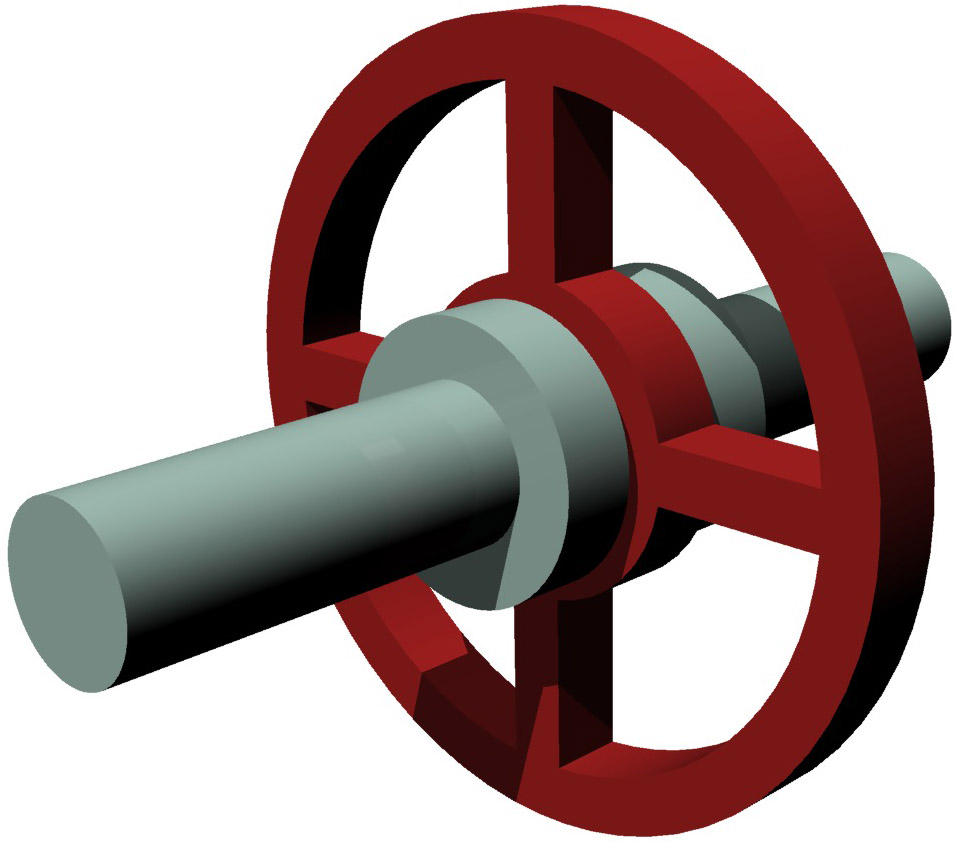}
\caption{A wheel on an axle.}
\label{Fig:Wheel}
\end{wrapfigure}

Next, a \emph{design} is a union of rigid bodies in space, with
disjoint interiors.
As an example of a design, we may take a wheel on an axle as shown in
\reffig{Wheel}.  We add two flanges to the axle to stop the wheel
from moving along the length of the axle.

A \emph{compound movement} $M$, of $D$, is a collection of movements,
one for each rigid body, which together keep the interiors of the
bodies disjoint.  In the example, spin the wheel in one direction
while spinning the axle in the other.
As this happens we may additionally translate both through space to
give a new compound movement.  In other words, if $M$ is a compound
movement of $D$ and $g$ is any movement then we define a new compound
movement $g \circ M$ by composing all of the movements of $M$ with
$g$, moment by moment.

Two compound movements $M$ and $N$ are \emph{equivalent} if there is a
movement $g$ so that $g \circ M = N$.  Thus, the two compound
movements of wheel plus axle discussed above are equivalent.

We now restrict the degrees of freedom of $D$ via the following axiom.
\begin{itemize}
\item[] \emph{Tracked}: There is a compound movement $M$ of $D$ so
  that every compound movement of $D$ is equivalent to a
  time-reparameterization of $M$.
\end{itemize}

The wheel plus axle satisfies tracked: we can rotate the wheel on the
axle and there is essentially no other way that they can move relative
to each other.  However, an axle with two wheels (again with flanges
holding them in place) does not satisfy tracked.
This axiom implies that in $D$, as in many geared machines, the bodies
may move relative to each other in an essentially unique fashion.  We
now distinguish the wheel and axle from the designs of Ferguson and
Van Deventer.
\begin{itemize}
\item[] \emph{Epicyclic}: For any bodies $A$ and $B$ in the design
  $D$, if $A$ and $B$ are in contact at some moment then the given
  movement of $B$, in the frame of reference of $A$, is not simple.
\end{itemize}
Fixing the axle pointwise, the wheel can only spin about its center,
which is a simple movement.  Hence the wheel plus axle is not
epicyclic.  Note that in ordinary mechanisms, the parts in contact
with the gear-box have simple movements, in the frame of reference of
the box.  Such mechanisms are not epicyclic.  In view of this axiom we
will refer to \emph{all} of the rigid bodies in our design as
\emph{gears}.

To simplify our search for a tracked, epicyclic three-component
design, we add the following axiom.
\begin{itemize}
\item[] \emph{Symmetry}: Any gear of $D$ may be taken to any other by
  a rigid motion preserving $D$.
\end{itemize}

To ensure that $D$ is tracked, we will make a series of assumptions,
at each stage reducing the possible compound movements of the design.
To ensure that $D$ is epicyclic we must do this without using any kind
of framework or gear-box, following Ferguson and Van Deventer.
Departing from their work, symmetry implies that all of our gears must
be identical.

\section{Topology}



In this section we prevent compound movements that move the gears
arbitrarily far from each other.  We first assume that each gear is a
circular ring.  Now, if no pair of rings is linked, then $D$ is
globally unlinked~\cite[Lemma~3.2]{freedman-skora}.  Thus we assume
that some pair of rings is linked.  For example, this rules out the
Borromean rings.
\begin{wrapfigure}[12]{l}{0.20\textwidth}
\centering
\vspace{-7pt}
\includegraphics[width=0.20\textwidth]{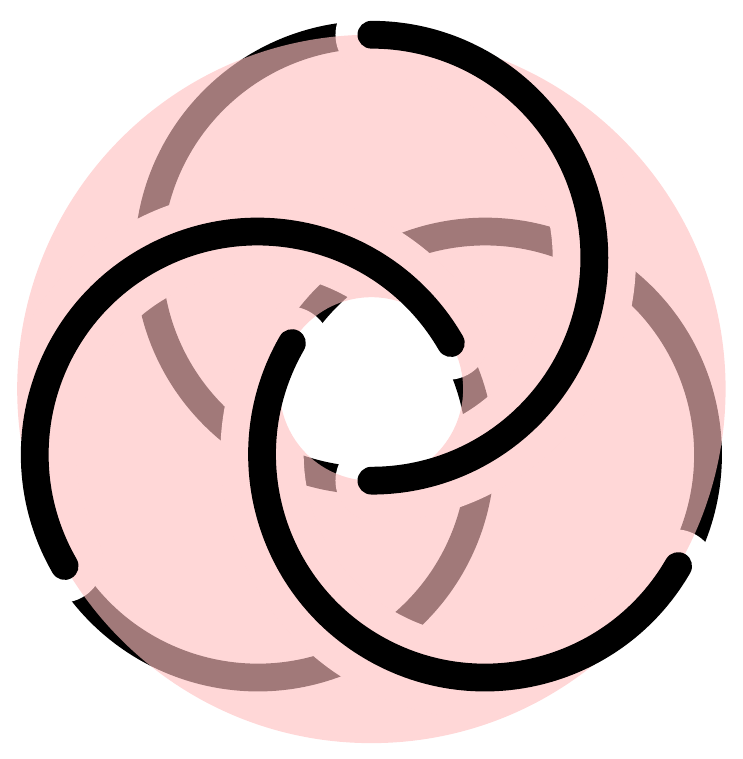}
\caption{The three component Hopf link, lying on a torus.} 
\label{Fig:Hopf}
\end{wrapfigure}
As suggested above, we further restrict our attention to designs with
only three rings.  The symmetry condition now implies that all three
pairs of rings are linked.  A geometric argument given in Appendix
\ref{Sec:HopfLink} proves the design $D$ is, topologically, the
three-component Hopf link.  Note that this rules out another natural
candidate: the minimally twisted three chain link~\cite{KaiserEtAl}.

Examining cases, the symmetry axiom implies that there is a rigid
motion, $g_0$, that acts via an order three permutation\footnote{That
  is, if we apply $g_0$ three times then we are back to where we
  started.}  on the components of $D$.  It follows that $g_0$ is a
rotation about a line in three-space of order three (i.e. by
$120^\circ$ in one or the other direction).
Our three-component Hopf link has just two symmetries of
order three.
Thus, we may assume that the rings of $D$ are equally spaced $(1,1)$
curves on a torus, as shown in \reffig{Hopf}.  

Consider the compound movement $M$ that rotates each of the rings
about its center, all with the same speed and handedness.  Note that
$M$ preserves the design and obeys the symmetry axiom.  Also, in the
frame of reference of any one ring, the movements of the other two are
not simple, so a design with this movement would be epicyclic.  In the
rest of the paper, we seek a shape for the individual rings that will
satisfy the tracked axiom, using essentially this compound movement
$M$.

\section{Maximal tori}
\label{Sec:Max}

\begin{wrapfigure}[15]{r}{0.60\textwidth}
\vspace{-25pt}
\centering 
\subfloat[Three rings arranged as a Hopf link, maximizing thickness
  among arrangements with 3--fold rotational symmetry.]
{
\includegraphics[width=0.29\textwidth]{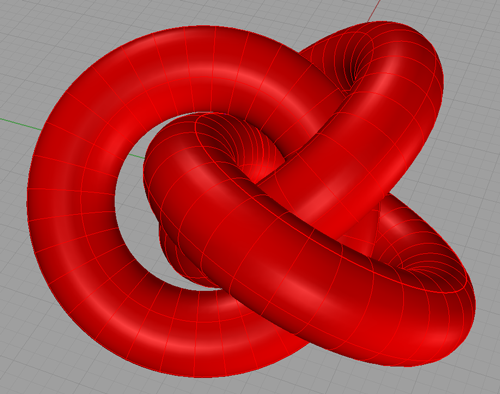}
\label{Fig:Max}
} \quad \subfloat[Variables for the optimization problem. 
  $U$ is perpendicular to $V$
  and lies in the $x,y$ plane. The line emanating from $C$ is parallel
  to the $z$ axis.]
{
\labellist
\small\hair 2pt
\pinlabel $C$ at 129 94
\pinlabel $U$ at 170 140
\pinlabel $V$ at 105 193
\pinlabel \textcolor{red}{$x$} at 253 112 
\pinlabel \textcolor{red}{$y$} at 8 200 
\pinlabel \textcolor{red}{$z$} at 83 243
\scriptsize
\pinlabel \textcolor{blue}{$r$} at 96 82
\pinlabel \textcolor{blue}{$\theta$} at 115 150
\tiny
\pinlabel \textcolor{blue}{$\phi$} at 183 123.5
\endlabellist

\includegraphics[width=0.26\textwidth]{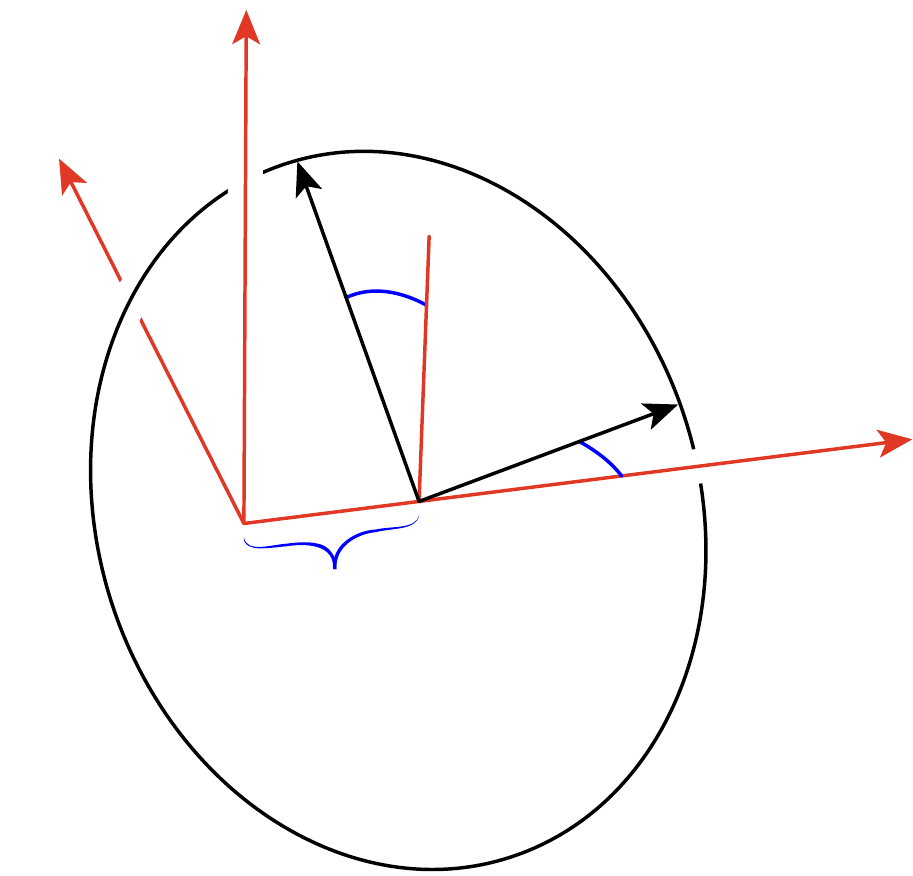}
\label{Fig:Optimization}
}
\caption{}
\label{Fig:OptimizationBig}
\end{wrapfigure}

In any design satisfying the tracked axiom, the gears must remain in
contact.  To arrange this we take three congruent rings, linked as in
\reffig{Hopf}, and gradually increase their thickness; think of
blowing up three tubular balloons.  Their orientations may change as
the rings bump each other, but we always maintain the $3$--fold
symmetry.  Eventually the rings will reach some maximum thickness,
giving an arrangement as in \reffig{Max}.

To find the positions of the rings at the maximal thickness we solve
an optimization problem.  The three circles, each of radius $1$, are
arranged so that rotation by $120^\circ$ around the $z$--axis permutes
them.  We assume that one of the circles has center on the $x$--axis,
at distance $r$ from the origin.  The plane containing this circle is
spanned by the vectors $U=(\cos(\phi), \sin(\phi),0)$ and
$V=(-\sin(\phi)\sin(\theta), \cos(\phi)\sin(\theta), \cos(\theta))$.
These parameters are displayed in \reffig{Optimization}.

Our aim is to maximize the distance between this circle and one of its
rotates, while preserving the topology.  Note that there is no closed
form for the minimal distance between two general circles~\cite{neff}.
Thus, a closed form solution for the optimal values of $r$, $\phi$,
and $\theta$ seems unlikely.  However, for our purposes a numerical
approximation is good enough.

Ian Agol~\cite{mathoverflow_question} points out that the optimal
configuration should have dihedral symmetry, and thus $\phi = 0$.
Here is the heuristic explanation.  First, at the maximal thickness,
each ring touches each of the other two in at least two points.  This
is because if there is only one point of contact, then the rings can
be rotated to remove that point of contact and allow further
thickening.  Second, two linked congruent tori that are tangent to
each other at two points should be related by a rotation of
$180^\circ$ that exchanges the two points.  This extra symmetry
property seems to hold (although we do not know a proof)
and implies the dihedral symmetry, which implies that $\phi = 0$.

We numerically obtain, for the other variables, values of $r =
0.4950197$ and $\theta = -0.8560281$.  At the maximum, the ring has
thickness of $0.3228837$ and the distance between two core circles is
twice this.  
We set up a coordinate system on the surface of each ring, a torus, as
follows.  Each torus is parameterized by $(\alpha, \beta) \in
(-\pi,\pi] \times (-\pi,\pi]$.  The parameter $\alpha$ is in the
    longitude direction, with $0$ at the vector $U$ and $V$ nearest in
    the positive direction.  The parameter $\beta$ is in the meridian
    direction, with $0$ at the biggest longitude (i.e. on the outside
    of the torus), and the direction $U \times V$ closest in the
    positive direction.  In these coordinates the points of contact of
    this torus with the others are at
\small
\[
(\alpha,\beta) = (-2.9419218,-1.2298655), \, (-1.9117269,2.9419218),  (1.9117269,-2.9419218), \, (2.9419218,1.2298655).
\]
\normalsize

\begin{wrapfigure}[12]{l}{0.58\textwidth}
\vspace{-24pt}
\centering 
\subfloat[Three rings arranged with 3--fold rotational symmetry.]
{
\includegraphics[width=0.222\textwidth]{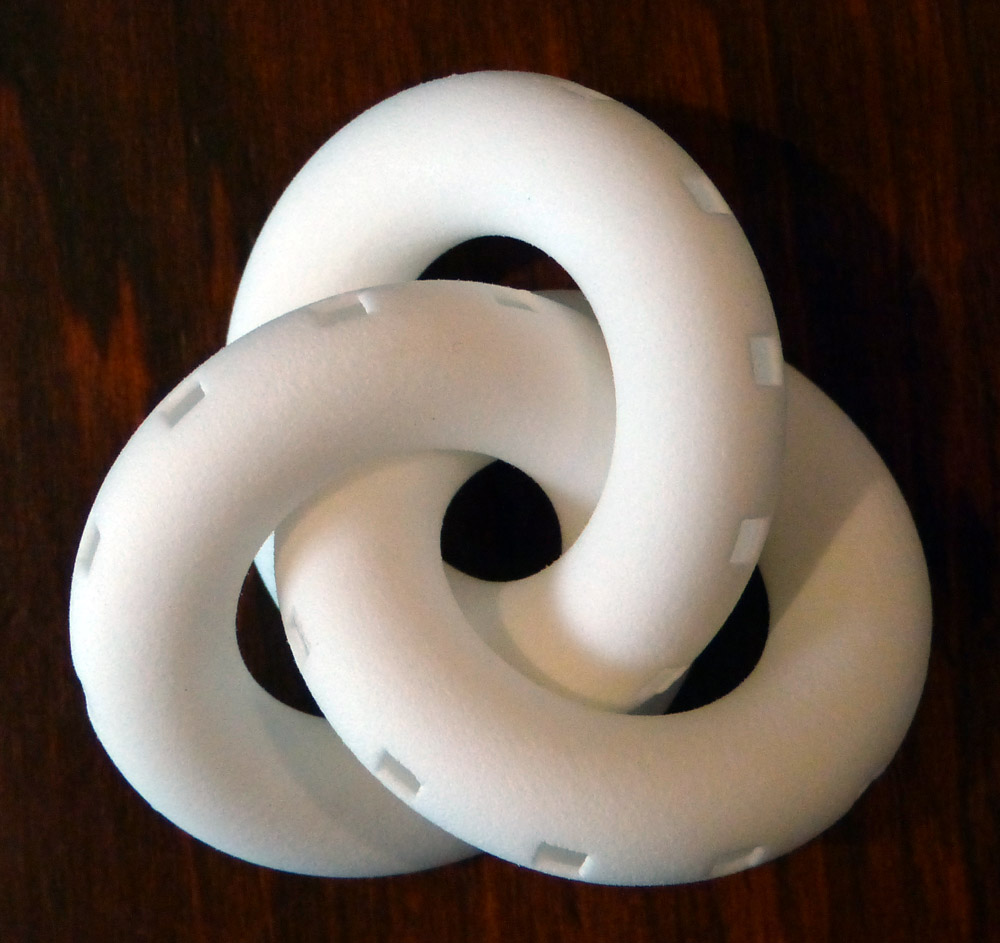}
\label{Fig:Symmetric}
}
\quad
\subfloat[If the three rings move out of the symmetrical position,
  then one or more of the six points of contact disappears.]  
{
\includegraphics[width=0.297\textwidth]{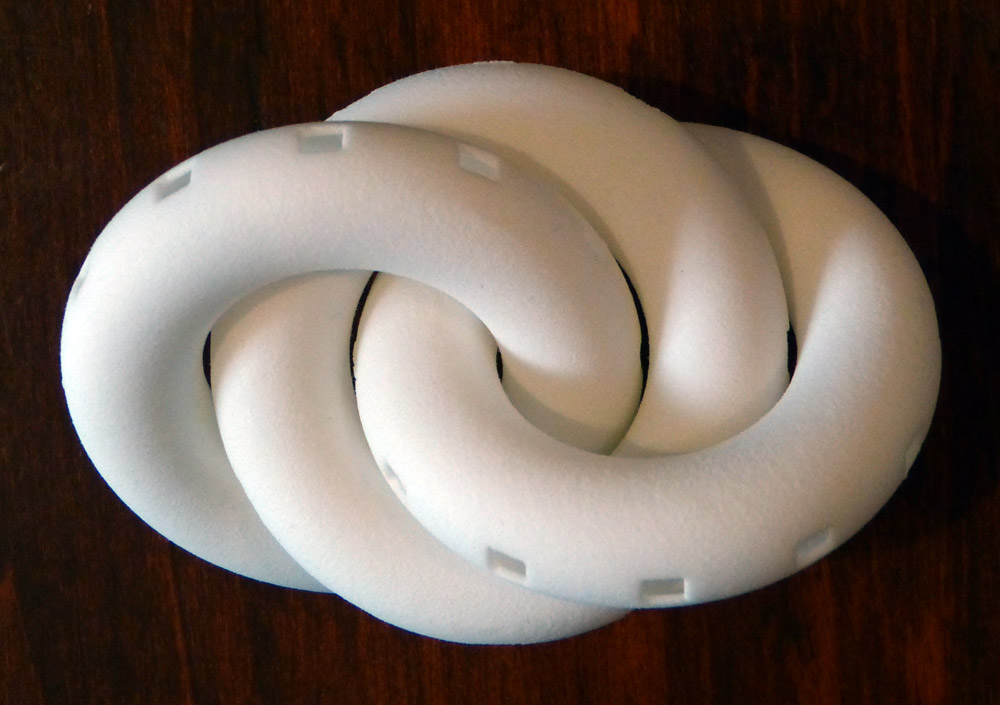}
\label{Fig:Asymmetric}
}
\caption{}
\label{Fig:printed tori}
\end{wrapfigure}

Note there are two points of contact with each of the other tori.  By
symmetry, the same is true for the other two tori as well.

Some of the symmetries in these four points come from the dihedral
symmetry of the configuration, but we currently do not know why the
number $2.9419218$ appears both as a value for $\alpha$ and for
$\beta$, nor why the sum of $1.2298655$ and $1.9117269$ appears to be
$\pi$.

If, at maximal thickness, the rings could move only along their core
circles, then the symmetry axiom would completely determine the
compound movement required by the tracked axiom.  Unfortunately this
is not the case: if we relax the $3$--fold symmetry then we can
further increase the thickness of the rings; see
Figures~\ref{Fig:Symmetric} and~\ref{Fig:Asymmetric}\footnote{We
  expect that maximizing the thickness of the rings in
  \reffig{Asymmetric} would restrict them to move only along their
  core circles.  An interesting future direction is to explore gearing
  in this configuration.}.  

\section{Cutting teeth}
\label{Sec:Carve}

We now add gear-teeth to the rings; this enforces the 3--fold rotational
symmetry, preventing the configuration of \reffig{Asymmetric}.  Recall
that symmetry requires our rings to all rotate at the same speed.

For unlinked gears with parallel axes, the tracked axiom (implied by
the fundamental law of gearing\footnote{The fundamental law of gearing
  states that the angular velocity ratio between two gears should
  remain constant as the gears rotate.}) can be obtained using
involutes of a circle for the shapes of the teeth; this was first
proposed by Euler~\cite{euler}.  The version of this problem for skew
axes is an area of active research.  We looked in particular at work
of Phillips~\cite{phillips_book, stachel} on spatial involute gearing.
We implemented his system, but ran into a serious problem.  For any
pair of axes, he gives a pair of surfaces that are tangent at the
contact point and curve away from the contact point towards the two
axes.  This nicely suits the case when the gear lies on the side of
the flank surface towards the axis, as in most applications.  However,
our gears link each other.  Thus, it can happen that the gear may lie
on the opposite side of Phillips' tooth flank.  It can even happen
that the flank surfaces are roughly cylindrical, tangent at a point,
curving in the same direction, and with different axes --- in this
case the two surfaces intersect in an `X' shape.

So we instead use an ad hoc method.  We use one gear to \emph{carve}
the shape of the mating gear.  In hindsight, our method is similar to
techniques going back at least to 1896~\cite[Figures~12--19]{Fellows}.
Here are the details.  We know the compound movement of the two
meshing gears.  We choose a point $p$ on the boundary of the first
gear.  We consider the path of $p$ in the frame of reference of the
second gear.  This path must be outside of, or on the boundary of, the
second gear; so we imagine this path (and many other similar paths)
carving the shape of the second gear, as if from a block of clay.

\begin{figure}[p]
\centering
\subfloat[Core circles of two of the tori, together with a set of gear
  teeth on the inside of the first ring.]
{
\includegraphics[width=0.398\textwidth]{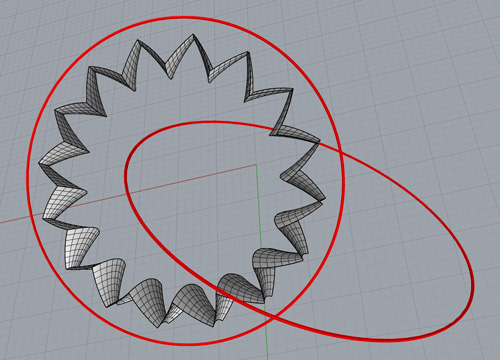}
\label{Fig:carve_gearing1.png}
}
\quad
\subfloat[A curve offset from a gear tooth, in order to leave a gap
  between the eventual gears.  Without this, the 3D printing process
  would produce a single merged object.]
{
\includegraphics[width=0.361\textwidth]{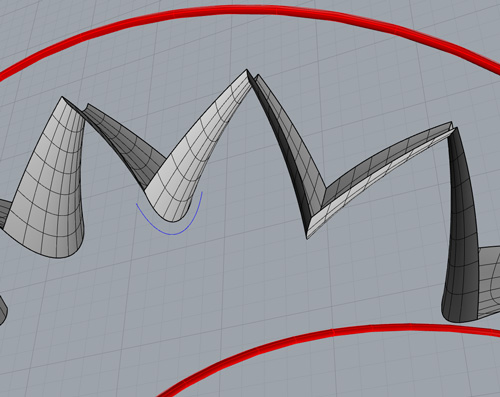}
\label{Fig:carve_gearing2.png}
}

\subfloat[Transforms of the curve in the frame of reference of the
  second ring.  Points on these curves closest to the core circle of
  the second ring are marked.]
{
\includegraphics[width=0.382\textwidth]{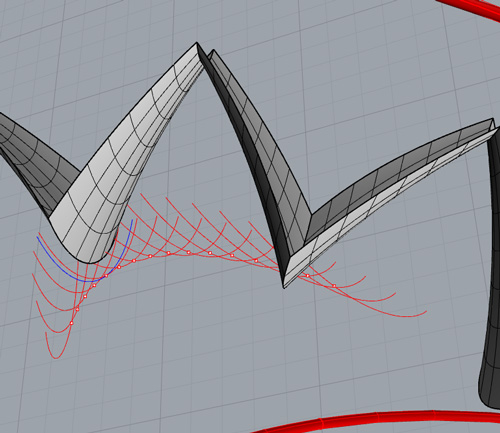}
\label{Fig:carve_gearing3.png}}
\quad
\subfloat[A curve fitted through these points.  By construction, this
  curve is on the boundary of the region carved out by the first
  gear.]
{
\includegraphics[width=0.382\textwidth]{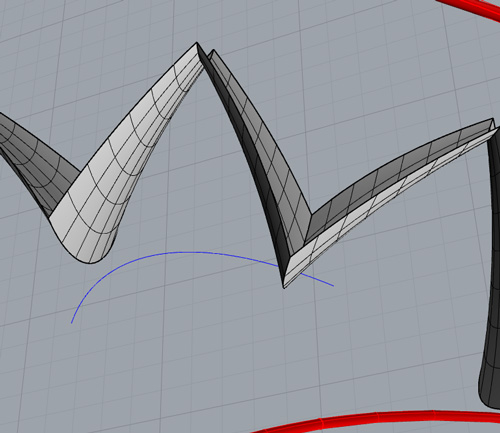}
\label{Fig:carve_gearing4.png}
}

\subfloat[Surfaces swept out along curves defined as in the previous
  steps give the tooth flanks for the second ring. The curve from
  \reffig{carve_gearing4.png} is the lower boundary of the left
  flank of the gear tooth.]
{
\includegraphics[width=0.403\textwidth]{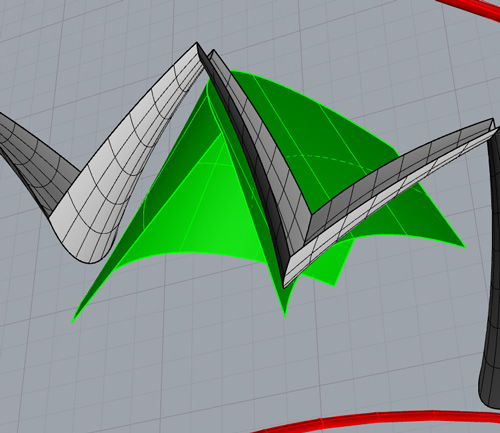}
\label{Fig:carve_gearing5.png}
}
\quad
\subfloat[The 3D printed triple gear.]
{
\includegraphics[width=0.361\textwidth]{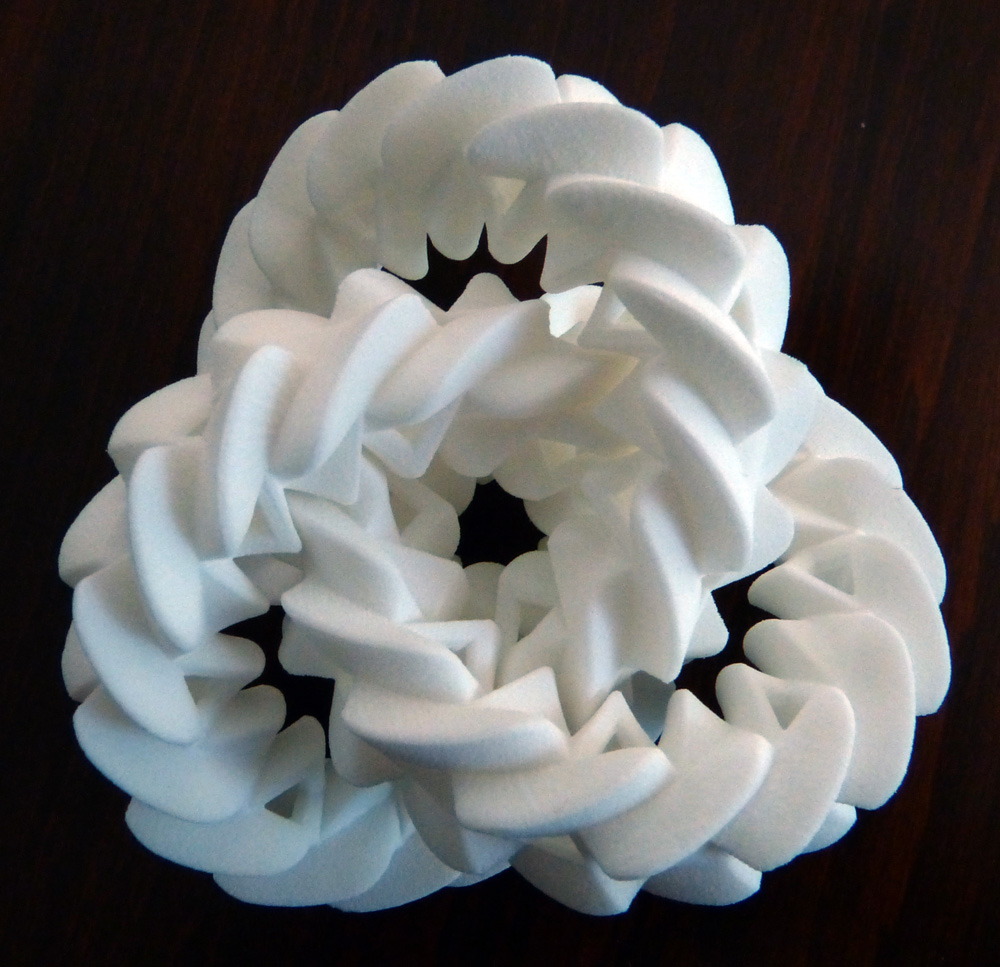}
\label{Fig:printed result}
}
\caption{All but the last of these diagrams are screenshots from
  Rhinoceros, the CAD software we used.  We also relied extensively on
  Rhinoceros' Python scripting interface.}
\label{Fig:carve gearing}
\end{figure}

Note that this process transforms a point of one gear into a excluded
path for the other gear; that is, the dimension increases by one.  If
we have a proposed $2$--parameter tooth flank for the first gear, then
this carves a $3$--parameter volume in the frame of reference of the
second gear.  The appropriate tooth flank for the second gear would
then be part of the boundary of that volume.  In practice, we start
with one set of gear teeth on the inside of one ring, and use it to
carve a set of gear teeth for a second.  See the steps outlined in
\reffig{carve gearing}.  The shapes of the first set of teeth are
given by linear functions in the toroidal coordinates; the exact
parameters were chosen with a fair amount of trial and error to work
well with the carving process.

Note that the $\beta$ values of the points of contact on the torus
determine longitudes of the torus that must have gearing on them. The
two values $\beta=-2.9419218$ and $\beta=2.9419218$ are very close to
each other near the smallest longitude on the inside of the ring.
These are so close that gear teeth from the other two rings intended
to mesh with one longitude will most likely interfere with the other.
For this reason, we merged these two longitudes into one gearing
pattern, as seen in \reffig{carve_gearing1.png}\footnote{The
  closeness of these two longitudes may cause problems in producing
  three gears that rotate at different speeds.}.  We follow the steps
shown in \reffig{carve gearing} to produce the tooth flanks.  As is
usual in gear design, we truncate to ensure a non-zero ``top land''
for the gear teeth.
We transport the tooth flanks using the various symmetries of the
design.  The flanks are joined to each other with connecting patches,
and the final closed surface is hollowed out; this allows recovery of
unused material from the 3D printing process.  

This completes our construction; \reffig{printed result} displays the
final result\footnote{See
  \url{http://www.youtube.com/watch?v=I9IBQVHFeQs} for a video of the
  triple gear in motion. It can be purchased at
  \url{http://shpws.me/mQ6F}. An STL file is available at
  \url{http://www.thingiverse.com/thing:66708}.}.  The gear teeth on
each ring are arranged in three groups, one along the inside longitude
of the ring and two along longitudes towards the outside of the
ring. There are the same number of teeth on each of these, so in
principle it might be possible to connect all of these together into a
single twisting gearing pattern wrapping around the torus. This would
have to be done very carefully, to avoid any unintended collisions
between the added parts and the existing gear teeth.

The 3-fold symmetric three-component Hopf link is one of many other
arrangements of linked circles that could lead to tracked, epicyclic
designs. In addition to the configuration shown in
\reffig{Asymmetric}, rotationally symmetric Hopf links with four or
five components should be possible. We have also considered a design
based on the minimally twisted six chain link~\cite{KaiserEtAl}, where
neighboring rings are geared while rings two steps around the chain
from each other are in contact but slide smoothly against each other.

\section{Driving the triple gear}

It is an interesting dexterity challenge to move the triple gear using
your hands.  Printed in PA 2200 nylon plastic by selective laser
sintering, friction wins over weight; if you hold two of the gears and
rotate them then the third is pushed out of the 3-fold symmetry
position until the mechanism seizes.  A finger positioned against the
third gear stops it moving out of place; then the design moves as
intended.

\begin{figure}[htbp]
\centering
\subfloat[The helical axle screws into the triple gear.]
{
\includegraphics[width=0.314\textwidth]{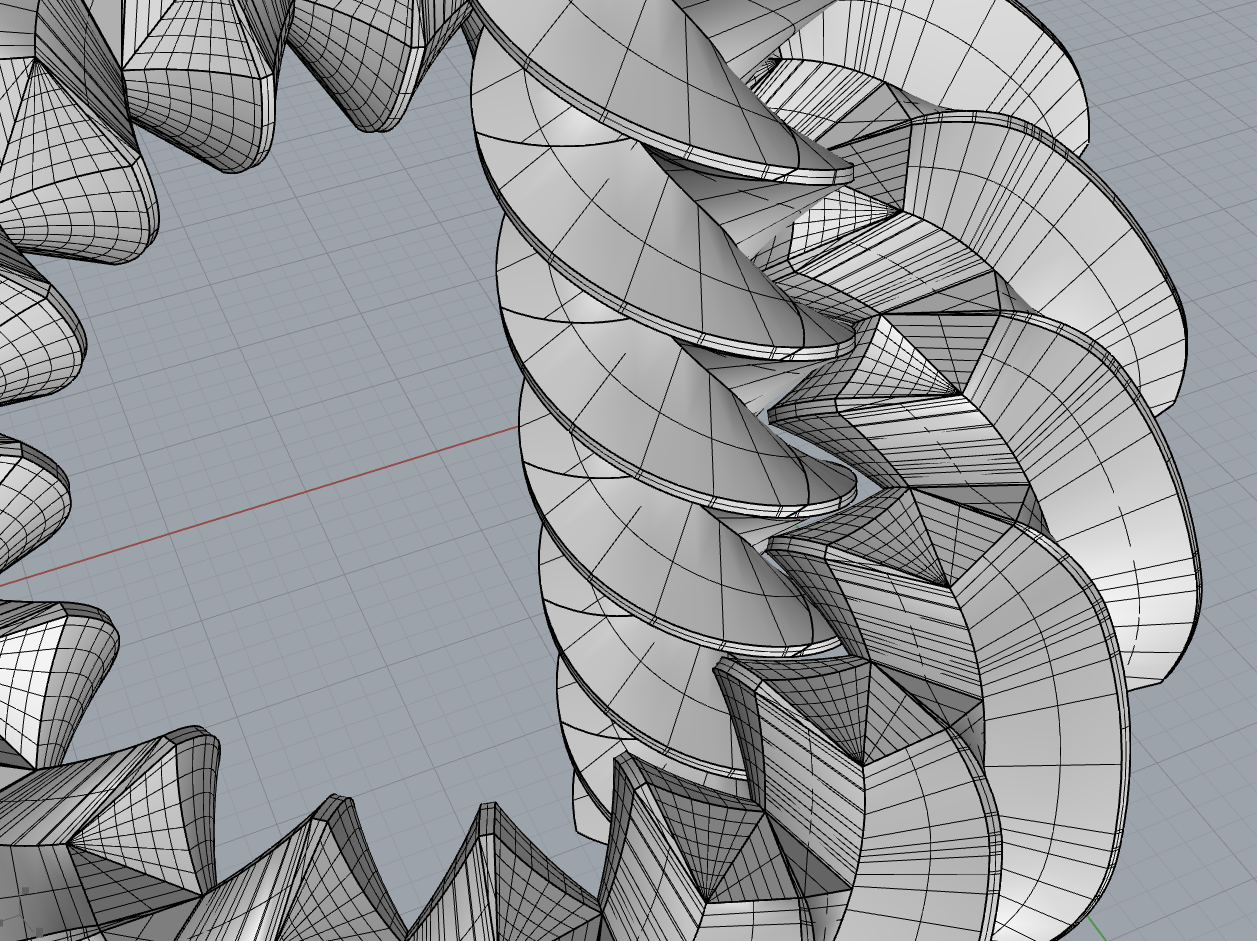}
\label{Fig:axle and one gear}
}
\quad
\subfloat[The baseplate and driving axle.]
{
\includegraphics[width=0.342\textwidth]{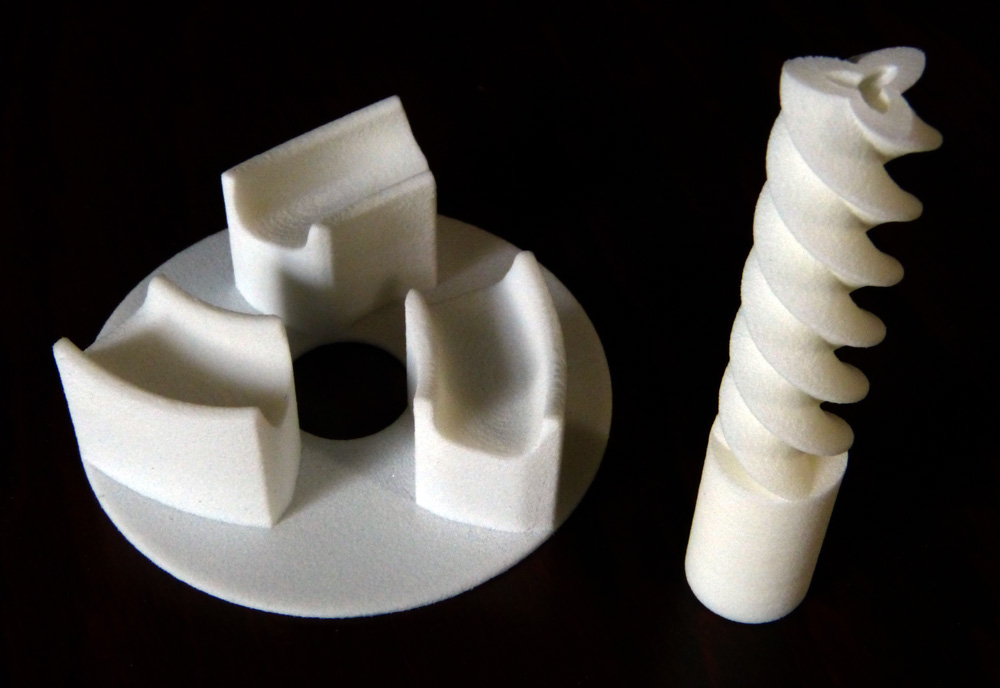}
\label{Fig:baseplate and axle}
}
\quad
\subfloat[The triple gear with baseplate and driving axle.]
{
\includegraphics[width=0.26\textwidth]{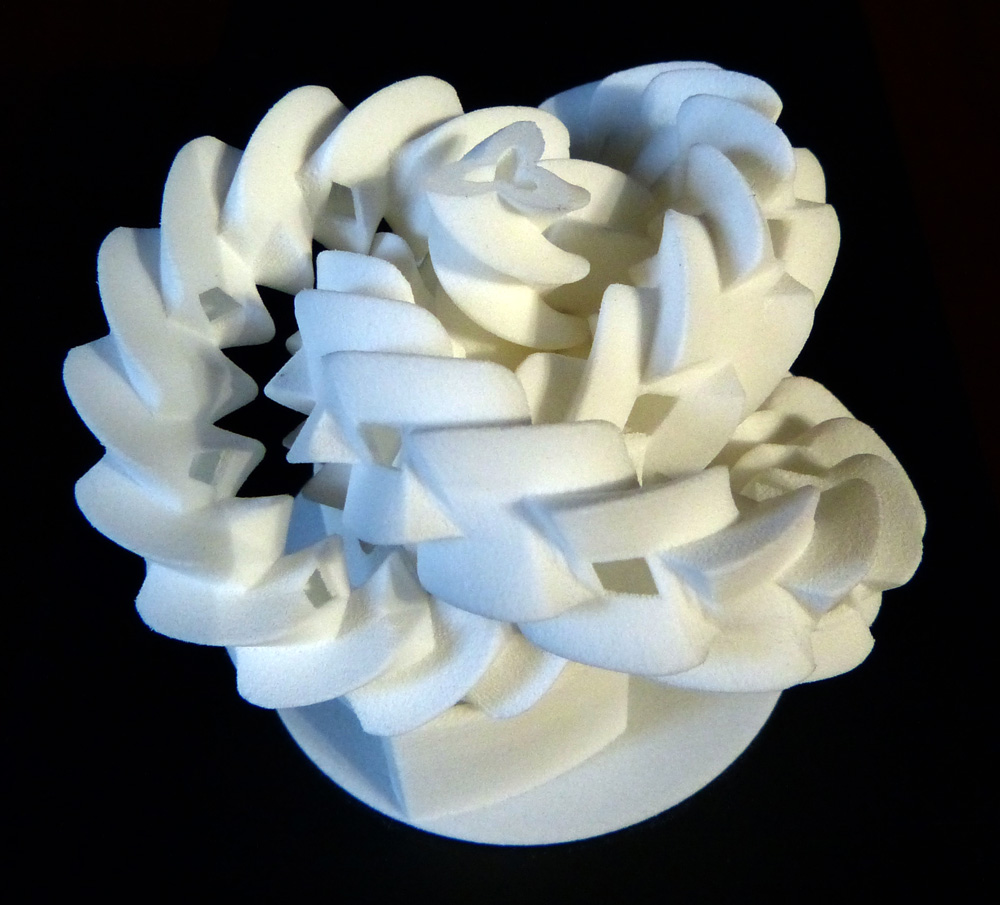}
}
\caption{}
\label{Fig:HelicalAxle}
\end{figure}

We have also designed a driving mechanism.  As shown in
\reffig{baseplate and axle}, it consists of a central helical axle,
powered by a motor, and a baseplate upon which the triple gear
rests\footnote{See \url{http://www.youtube.com/watch?v=QhXjevOY_uk}
  for a video of the powered triple gear. We thank Adrian Goldwaser
  for initial prototyping and Stuart Young for further prototyping and
  construction of the motorized base. The baseplate and axle can be
  purchased at \url{http://shpws.me/nIjE}.  An STL file is available
  at \url{http://www.thingiverse.com/thing:78218}.}.  The tooth
cross-section for the helical axle is designed using essentially the
same carving technique as in \refsec{Carve}.  Since the helix is
invariant under a screw motion it can be inserted into the the triple
gear without having to move the gears; it need not be printed in
place.

The axle, with a rotational movement, can be added to the design $D$
with its original movement; this combination is again
epicyclic.  However, in reality this new design tends to prefer the
compound movement where all of the circular gears are at rest with
respect to the axle: that is, everything rotates about the axle at the
same speed.  This is definitely not epicyclic!  The three cradles on
the baseplate prevent this trivial compound motion.

Since rotational and translational movement of the helical axle are
equivalent, we can also drive the triple gear by translating the axle
through it.  Thus we can place the triple gear at the top of a long
axle and let it fall downwards under gravity.  However, without the
baseplate the triple gear only sometimes moves in the desired
epicyclic fashion.  It should be possible to produce a very large
circular axle by gently bending the long axle into a ring.  Holding
the circular axle vertically and, pulling the axle hand-over-hand, the
triple gear falls endlessly down the other side.


\setlength{\baselineskip}{13pt} 

\bibliographystyle{hamsplain}
\bibliography{triplegearbib}

\appendix
\section{Postscript on paradoxes}

After submission, but before publication, we stumbled across the work
of Jacques Maurel\footnote{See
  \url{http://www.jacquesmaurel.com/gears}.  His video
  \url{http://www.youtube.com/watch?v=AwzL7Z_50zc} inspired
  this appendix.}.  He explains, via many beautiful videos, Lindsay's
\emph{paradoxical gears}~\cite{Lindsay}.  We adapt these to find three
gears with parallel axes, all in contact and all rotating in the same
direction.  However, these gears are not in fact planar; thus they
narrowly avoid contradicting the second sentence of our abstract.


\begin{wrapfigure}[24]{l}{0.20\textwidth}
\vspace{-5pt}
\caption{Horizontal slices of the paradoxical parallel axis gears.}
\centering
\vspace{-5pt}
\includegraphics[width=0.20\textwidth]{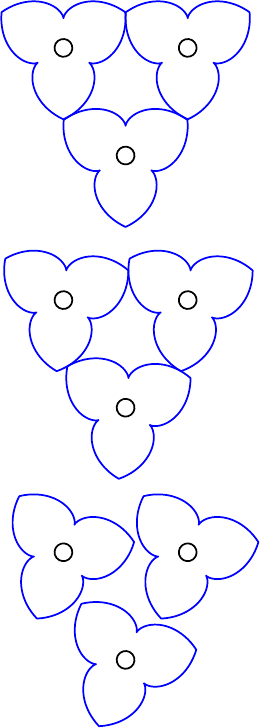}
\label{Fig:Paradoxical}
\end{wrapfigure}

The three gear are helical screws, similar to the driving axle shown
in \reffig{axle and one gear}.  Their axes meet each horizontal plane
at the vertices of an equilateral triangle.  Each horizontal
cross-section of each screw is an involute planar gear, as discussed
in \refsec{Carve}, but with only three teeth.  
See \reffig{Paradoxical} for a sequence of cross-sectional views.
Descending, the cross-sections rotate, all in the same direction.
Alternatively we can view the three pictures as an animation of a
cross-section at a fixed height.


Note that in ordinary planar gearing, teeth of meshing gears meet
close to the \emph{center line} $L$: the straight line between the
gear centers.  The normal to the teeth at the point of contact is
roughly perpendicular to $L$.  Thus, using the cross-product to
compute the torque on the follower, we find that the driver and
follower rotate in opposite directions.

Consider the first cross-section of the paradoxical gears, shown in
\reffig{Paradoxical}.  Again there is a triple of involute gears, in
point contact.
However, the point of contact of the upper two gears is much further
from their center line, $L$.  At the point of contact the common
normal to the teeth is parallel to $L$.  Thus the follower rotates in
the same direction as the driver.

By the time we get to the third picture in the animation, the gears
are no longer in contact.  However, there are horizontal slices of the
screws at other heights which have come into contact, reproducing the
first picture.  So the helical screws, as a whole, remain in contact
and produce the desired paradoxical behavior.


\section{The three-component Hopf link}
\label{Sec:HopfLink}

\begin{no_num_thm} 
Let $L$ be a three-component link where every component is a round
circle and every pair of circles is linked.  Then $L$ is the
three-component Hopf link.
\end{no_num_thm}

\begin{proof}
We work in $S^3 = \{(z,w) \in \CC^2 : |z|^2 + |w|^2 = 1\}$.  Let $Z =
\{(z,0) : |z| = 1\}$ and $W = \{(0,w) : |w| = 1\}$ be the $z$-- and
$w$--axes.  After a M\"{o}bius transformation, we may assume that the
first and second components of the link $L$ are the axes $Z$ and $W$.
Let $L_3$ be the third component of $L$ and note that this is again a
round circle.

Consider the set of great spheres that contain $Z$.  Each such sphere
is divided into two hemispheres by $Z$.  Each sphere meets $W$ in two
points, one in each hemisphere.  So the set of hemispheres $\{H_w :
|w| = 1\}$ is indexed by the points of $W$.

Since $L_3$ is round, and since it links $Z$, we deduce that $L_3$ intersects
each great sphere containing $Z$ in precisely two points, once in each
hemisphere.  Thus, as we go around $L_3$, we hit every hemisphere
$H_w$ in cyclic order.  Similarly, as we go around $L_3$, we hit the
hemispheres $H_z$ in cyclic order.  Note that $H_z \cap H_w$ is the
geodesic segment from $(z,0)$ to $(0,w)$, parameterized by $(z \cdot
\cos(\theta), w \cdot \sin(\theta))$ for $\theta \in [0,\pi/2]$.  We
can deform $L_3$ by sliding each point along its respective geodesic
segment until it is at $\frac{\sqrt{2}}{2}(z, w)$.  Thus $L_3$ becomes
the $(1,1)$ curve on the torus $|w| = |z|$.  So $Z$, $W$ and $L_3$ are
three fibers of the Hopf fibration of $S^3$, and our link is the
three-component Hopf link.
\end{proof}

\end{document}